\documentclass[11pt]{amsart}

\pdfoutput=1

\usepackage{amsmath, amsthm, amssymb, graphicx, url}

\newcommand\kh{\textit{Kh}}
\newcommand\kho{\textit{Kh}^\prime}
\newcommand\khr{\overline{\textit{Kh}}}
\newcommand\khor{\overline{\textit{Kh}}^\prime}

\newcommand\zz{\mathbb{Z}}

\newcommand\zzt{\zz/2\zz}

\newcommand\cdr{\overline{C}(\dd)}
\newcommand\cdrr{\widetilde{C}(\dd)}
\newcommand\cdrrr{\widehat{C}(\dd)}
\newcommand\dd{\mathcal{D}}
\newcommand\ee{\mathcal{E}}
\newcommand\vv{\mathcal{V}}
\newcommand\vo{\overline{V}}
\newcommand\zx{\mathbb{Z}\langle x_1, \dots, x_n \rangle}
\newcommand\vt{\widetilde{V}}
\newcommand\vh{\widehat{V}}

\newtheorem{theorem}{Theorem}
\newtheorem{cor}[theorem]{Corollary}
\newtheorem{lem}[theorem]{Lemma}
\newtheorem{prop}[theorem]{Proposition}

\theoremstyle{remark}
\newtheorem*{rem}{Remark}

\title[Odd Khovanov homology is mutation invariant]{Odd Khovanov homology is mutation invariant}

\author[Jonathan M.\ Bloom]{Jonathan M.\ Bloom}
\address{Department of Mathematics,
Columbia University,
New York, NY 10027, USA}
\email{jbloom@math.columbia.edu}

\begin{document}

\begin{abstract}
We define a link homology theory that is readily seen to be both isomorphic to reduced odd Khovanov homology and fully determined by data impervious to Conway mutation.  This gives an elementary proof that odd Khovanov homology is mutation invariant over $\zz$, and therefore that Khovanov homology is mutation invariant over $\zzt$.  We also establish mutation invariance for the entire Ozsv\'ath-Szab\'o spectral sequence from reduced Khovanov homology to the Heegaard Floer homology of the branched double cover.
\end{abstract} 

\maketitle

\section{Introduction}

To an oriented link $L \subset S^3$, Khovanov \cite{kh1} associates a bigraded homology group $\kh(L)$ whose graded Euler characteristic is the unnormalized Jones polynomial.  This invariant also has a reduced version $\khr(L,K)$, which depends on a choice of marked component $K$.  While the Jones polynomial itself is insensitive to Conway mutation, Khovanov homology generally detects mutations that swap arcs between link components \cite{we}.  Whether the theory is invariant under component-preserving mutation, and in particular for knots, remains an interesting open question, explored in \cite{bar} and \cite{kof}.  No counterexamples exist with $\leq 14$ crossings, although Khovanov homology does distinguish knots related by genus 2 mutation \cite{dun}, whereas the (colored) Jones polynomial does not.

In 2003, Ozsv\'ath and Szab\'o \cite{osz12} introduced a link surgery spectral sequence whose $E^2$ term is $\khr(L; \zzt)$ and which converges to the Heegaard Floer homology of the branched double cover $\Sigma(L)$.  In search of a candidate for the $E^2$ page over $\zz$, Ozsv\'ath, Rasmussen and Szab\'o \cite{orsz} developed odd Khovanov homology $\kho(L)$, a theory whose mod 2 reduction coincides with that of Khovanov homology.  While the reduced version $\khor(L)$ categorifies the Jones polynomial as well, it is independent of the choice of marked component and determines $\kho(L)$ according to the equation
\begin{align}
\label{eqn:red}
\kho_{m,s}(L) \cong \khor_{m,s-1}(L) \oplus \khor_{m,s+1}(L).
\end{align}
In contrast to Khovanov homology for links, we prove:
\begin{theorem}
\label{thm:main}
Odd Khovanov homology is mutation invariant.  Indeed, connected mutant link diagrams give rise to isomorphic odd Khovanov complexes.
\end{theorem}
\begin{cor}
\label{cor:1}
Khovanov homology over $\zzt$ is mutation invariant.
\end{cor}
\noindent We expect similar results to hold for genus 2 mutation.  Wehrli has announced an independent proof of Corollary \ref{cor:1} for component-preserving mutations, using an approach outlined by Bar-Natan.

Mutant links in $S^3$ have homeomorphic branched double covers.  It follows that the $E^\infty$ page of the link surgery spectral sequence is also mutation invariant.  Baldwin \cite{b} has shown that all pages $E^i$ with $i \geq 2$ are link invariants, to which we add:
\begin{theorem}
\label{thm:2}
The $E^i$ page of the Ozsv\'ath-Szab\'o spectral sequence is mutation invariant for $i \geq 2$.  Indeed, connected mutant link diagrams give rise to isomorphic filtered complexes.
\end{theorem}
\noindent
Note that Khovanov homology, even over $\zzt$, is not an invariant of the branched double cover itself \cite{w}.  We arrived at the above results in the course of work on the relationship between odd Khovanov homology and the monopole Floer homology of branched double covers \cite{jmb1}.

\subsection{Organization}  In Section \ref{sec:reform}, we give a new construction, which associates a bigraded complex $(\cdrr,\tilde\partial_\epsilon)$ to a connected, decorated link diagram $\dd$.  The construction takes as input a particularly small set of mutation invariant data derived from $\dd$, as summarized in Proposition \ref{prop:thrift}.  Nevertheless, in Section \ref{sec:oddkh} we prove:
\begin{prop}
\label{prop:id}
$(\cdrr, \tilde\partial_\epsilon)$ is canonically isomorphic to the reduced odd Khovanov complex $(\cdr, \bar \partial_\epsilon)$.
\end{prop}
\noindent This establishes Theorem \ref{thm:main} and verifies that our construction yields a well-defined link invariant.

In Section \ref{sec:bdc}, we consider a surgery diagram for the branched double cover $\Sigma(L)$ given by a framed link $\mathbb{L} \subset S^3$ with one component for each crossing of $\dd$.  Ozsv\'ath and Szab\'o associate a filtered complex to $\dd$ by applying the Heegaard Floer functor to a hypercube of 3-manifolds and cobordisms associated to various surgeries on $\mathbb{L}$ (see \cite{osz12}).  The link surgery spectral sequence is then induced by standard homological algebra.  In Proposition \ref{prop:surg}, we show that $\mathbb{L}$ itself is determined up to framed isotopy by the same mutation invariant data as $(\cdrr, \tilde\partial_\epsilon)$, establishing Theorem \ref{thm:2}.  We conclude with a remark on the original motivation for the construction of $(\cdrr, \tilde\partial_\epsilon)$ coming from branched double covers.  From this perspective, Theorems \ref{thm:main} and \ref{thm:2} are both immediate corollaries of the fact that taking the branched double cover of a link destroys all evidence of mutation.

\subsection{Acknowledgements}
I wish to thank John Baldwin, Josh Greene, Adam Levine, Ina Petkova, and my advisor Peter Ozsv\'ath for insightful discussions.  I was especially influenced by Greene's construction of Heegaard diagrams for branched double covers \cite{gr} and Baldwin and Plamenevskaya's construction of the reduced Khovanov homology of open book decompositions \cite{b3}.  I recently learned that Baldwin had programmed an algorithm \cite{b2} to compute the latter that is similar to the mod 2 reduction of the following construction.

\section{A thriftier construction of reduced odd Khovanov homology}
\label{sec:reform}

\begin{figure}[b]
  \centering
  \includegraphics{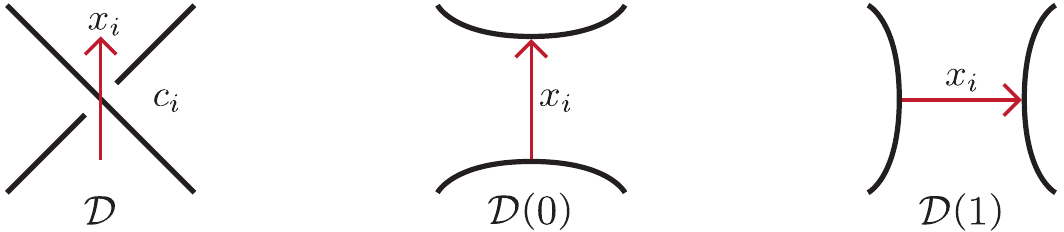}
  \caption{\textbf{Oriented resolution conventions.} The arrow $x_i$ at crossing $c_i$ remains fixed in a 0-resolution and rotates $90^\circ$ clockwise in a 1-resolution.  To see the other choice for the arrow at $c_i$, rotate the page $180^\circ$.  Mutation invariance of the Jones polynomial follows from the rotational symmetries of these tangles in $\mathbb{R}^3$.}
\label{fig:res_conv}
\end{figure}

Given an oriented link $L$, fix a connected, oriented link diagram $\dd$ with crossings $c_1, \dots, c_n$.  Let $n_+$ and $n_-$ be the number of positive and negative crossings, respectively.  We use $\vv(\dd)$ and $\ee(\dd)$ to denote the sets of vertices and edges, respectively, of the hypercube $\{0,1\}^n$, with edges oriented in the direction of increasing weight.  Decorate each crossing $c_i$ with an arrow $x_i$ which may point in one of two opposing and parallel directions.  The arrows $x_1, \dots, x_n$ appear in each complete resolution $\dd(I)$ as oriented arcs between circles according to the conventions in Figure \ref{fig:res_conv}.

Recall that a planar link diagram $\dd$ admits a checkerboard coloring with white exterior, as illustrated in Figure \ref{fig:mutation}.  The black graph $\mathcal{B}(\dd)$ is formed by placing a vertex in each black region and joining two vertices whenever there is a crossing in $\dd$ that is incident to both regions.  Given a connected spanning tree $\mathcal{T} \subset \mathcal{B}(\dd)$, we can form a resolution of $\dd$ consisting of only one circle by merging precisely those black regions which are incident along $\mathcal{T}$.  In particular, all connected diagrams admit at least one such resolution.

With these preliminaries in place, we now give a recipe for associating a bigraded chain complex $(\cdrr, \tilde\partial_\epsilon)$ to the decorated diagram $\dd$.  The key idea is simple.  Think of each resolution of $\dd$ as a connected, directed graph whose vertices are the circles and whose edges are the oriented arcs.  While the circles merge and split from one resolution to the next, the arcs are canonically identified throughout.  So we use the exact sequence
\begin{align*}
\{ \text{cycles}\} \hookrightarrow \zz\langle\text{arcs}\rangle \xrightarrow{d} \zz\langle\text{circles}\rangle \xrightarrow{} \zz \xrightarrow{} 0
\end{align*}
of free Abelian groups to suppress the circles entirely and instead keep track of the cycles in each resolution, thought of as relations between the arcs themselves.

We begin the construction by fixing a vertex $I^* = (m_1^*, \dots, m_n^*) \in \vv(\dd)$ such that the resolution $\dd(I^*)$ consists of only one circle $S$.  To each pair of oriented arcs $(x_i, x_j)$ in $\dd(I^*)$ we associate a linking number $a_{ij} \in \{0, \pm 1\}$ according to the symmetric convention in Figure 2.  We set $a_{ii} = 0$.  For each $I \in \vv(\dd)$, we have an Abelian group 
$$\vt(\dd(I)) = \zz \langle x_1,\dots, x_n \, | \, r^I_1, \dots, r^I_n \rangle$$
presented by relations
\begin{align}
\label{eqn:rel}
r^I_i &= \left\{
    \begin{array}{l}
        x_i + \displaystyle{\sum_{\{j \, | \, m_j \neq m_j^*\}}} (-1)^{m^*_j} a_{ij} x_j \quad \text{if} \quad m_i = m_i^* \\
        \hspace{8.2mm} \displaystyle{\sum_{\{j \, | \, m_j \neq m^*_j\}}} (-1)^{m_j^*} a_{ij} x_j \quad \text{if} \quad m_i \neq m_i^*.
    \end{array} \right.
\end{align}
Indeed, these relations generate the cycles in the graph of circles and arcs at $\dd(I)$ (see Lemma \ref{lem:1}).
\begin{figure}[b]
  \centering
  \includegraphics{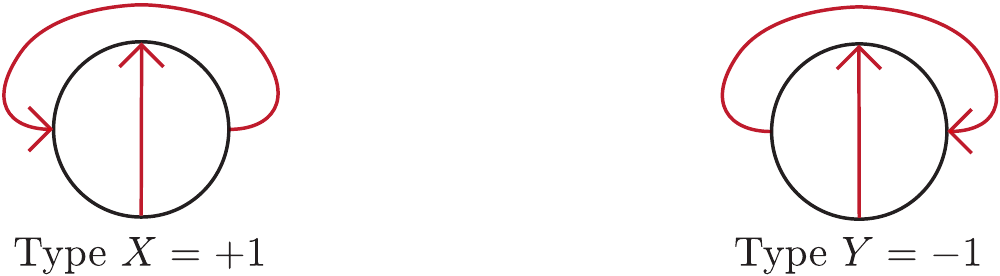}
  \caption{\textbf{Linking number conventions.} Two arcs are linked in $\dd(I^*)$ if their endpoints are interleaved on the circle.  Any linked configuration is isotopic to one of the above on the 2-sphere $\mathbb{R}^2 \cup \{\infty\}$.}
\label{fig:link_conv}
\end{figure}

To an edge $e \in \ee(\dd)$ from $I$ to $J$ given by an increase in resolution at $c_i$, we associate the map
$$\tilde\partial^I_J : \Lambda^* \vt(\dd(I)) \to \Lambda^* \vt(\dd(J))$$ defined by
\begin{align*}
\tilde\partial^I_J (u) &= \left\{
    \begin{array}{l}
        x_i \wedge u \quad \text{if} \quad x_i = 0 \in \vt(\dd(I)) \\
        \hspace{7.63mm} u \quad \text{if} \quad x_i \neq 0 \in \vt(\dd(I)). \\ 
    \end{array} \right.
\end{align*}
Extending by zero, we may view each of these maps as an endomorphism of the group
$$\cdrr = \bigoplus_{I \in \vv(\dd)} \Lambda^* \vt(\dd(I)).$$
Consider a 2-dimensional face from $I$ to $J$ in the hypercube corresponding to an increase in resolution at $c_i$ and $c_j$.  The two associated composite maps in $\cdrr$ commute up to sign and vanish identically if and only if the arcs $x_i$ and $x_j$ in $\dd(I)$ are in one of the two configurations in Figure \ref{fig:link_conv}, denoted type $X$ and type $Y$.  Note that we can distinguish a type X face from a type Y face without reference to the diagram by checking which of the relations $x_i \pm x_j = 0$ holds at each of the two vertices in $\cdrr$ strictly between $I$ and $J$.

A Type $X$ edge assignment on $\cdrr$ is a map $\epsilon : \ee(\dd) \to \{\pm 1\}$ such that the product of signs around a face of Type $X$ or Type $Y$ agrees with the sign of the linking convention in Figure \ref{fig:link_conv} and such that, after multiplication by $\epsilon$, every face of $\cdrr$ anticommutes.    Such an assignment defines a differential $\tilde\partial_\epsilon : \cdrr \to \cdrr$ by
\begin{align*}
\tilde\partial_\epsilon(v) = \sum_{\{e \in \ee(\dd), J \in \vv(\dd) \, | \, \text{$e$ goes from $I$ to $J$}\}} \epsilon(e) \cdot \tilde\partial^I_J(v)
\end{align*}
for $v \in \Lambda^* \vt(\dd(I))$.  Type $X$ edge assignments always exist and any two yield isomorphic complexes, as do any two choices for the initial arrows on $\dd$.  We can equip $(\cdrr, \tilde\partial_\epsilon)$ with a bigrading that descends to homology and is initialized using $n_\pm$ just as in \cite{orsz}.  The bigraded group $\cdrr$ and maps $\tilde\partial^I_J$ are constructed entirely from the numbers $a_{ij}$, $m_i^*$, and $n_\pm$. Thus, up to isomorphism:

\begin{prop}
\label{prop:thrift}
The bigraded complex $(\cdrr, \tilde\partial_\epsilon)$ is determined by the linking matrix of any one-circle resolution of $\dd$, the vertex of this resolution, and the number of positive and negative crossings.
\end{prop}

\begin{proof}[Proof of Theorem 1]
The following argument is illustrated in Figure \ref{fig:mutation}.  Given oriented, mutant links $L_1$ and $L_2$, fix a corresponding pair of oriented, connected diagrams $\dd$ and $\dd^\prime$ for which there is a circle $C$ exhibiting the mutation.  This circle crosses exactly two black regions of $\dd$, which we connect by a path $\Gamma$ in $\mathcal{B}(\dd)$.  To simplify the exposition, we will assume there is a crossing between the two strands of $\dd$ in $C$, so that $\Gamma$ may be chosen in $C$.  Extend $\Gamma$ to a spanning tree $\mathcal{T}$ to obtain a resolution $\dd(I^*)$ with one circle.  The natural pairing of the crossings of $\dd$ and $\dd^\prime$ induces an identification $\vv(\dd) \cong \vv(\dd^\prime)$.  The resolution $\dd^\prime(I^*)$ corresponds to a mutation of $\dd(I^*)$ and also consists of one circle $S$.  We can partition the set of arcs inside $C$ into those which go across $S$ (in dark blue) and those which have both endpoints on the same side of $S$ (in red).  Since this division is preserved by mutation, the mod 2 linking matrix is preserved as well.

Arrows at the crossings of $\dd$ orient the arcs of $\dd(I^*)$, which in turn orient the arcs of $\dd^\prime(I^*)$ via the mutation.  To preserve the linking matrix at $I^*$ with sign, we modify the arcs of $\dd^\prime(I^*)$ as follows.  Let $A = \{\text{arcs in } C \text{ and in } S \}$ and $B = \{\text{arcs in } C \text{ and not in } S \}$.  We reverse those arcs of $\dd^\prime(I^*)$ that lie in $A$, $B$, or $A \cup B$, according to whether the mutation is about the $z$-, $y$-, or $x$-axis, respectively (as represented at right in Figure \ref{fig:mutation}).  We then select a corresponding set of arrows on $\dd^\prime$.  Note that it may also be necessary to switch the orientations of both strands inside $C$ so that $\dd^\prime$ will be consistently oriented.  In any case, the number of positive and negative crossings is unchanged.  Propositions \ref{prop:id} and \ref{prop:thrift} now imply the theorem for $(\cdr, \bar\partial_\epsilon)$.  The unreduced odd Khovanov complex is determined by the reduced complex, just as in \eqref{eqn:red}.
\end{proof}

\begin{figure}
  \centering
  \includegraphics{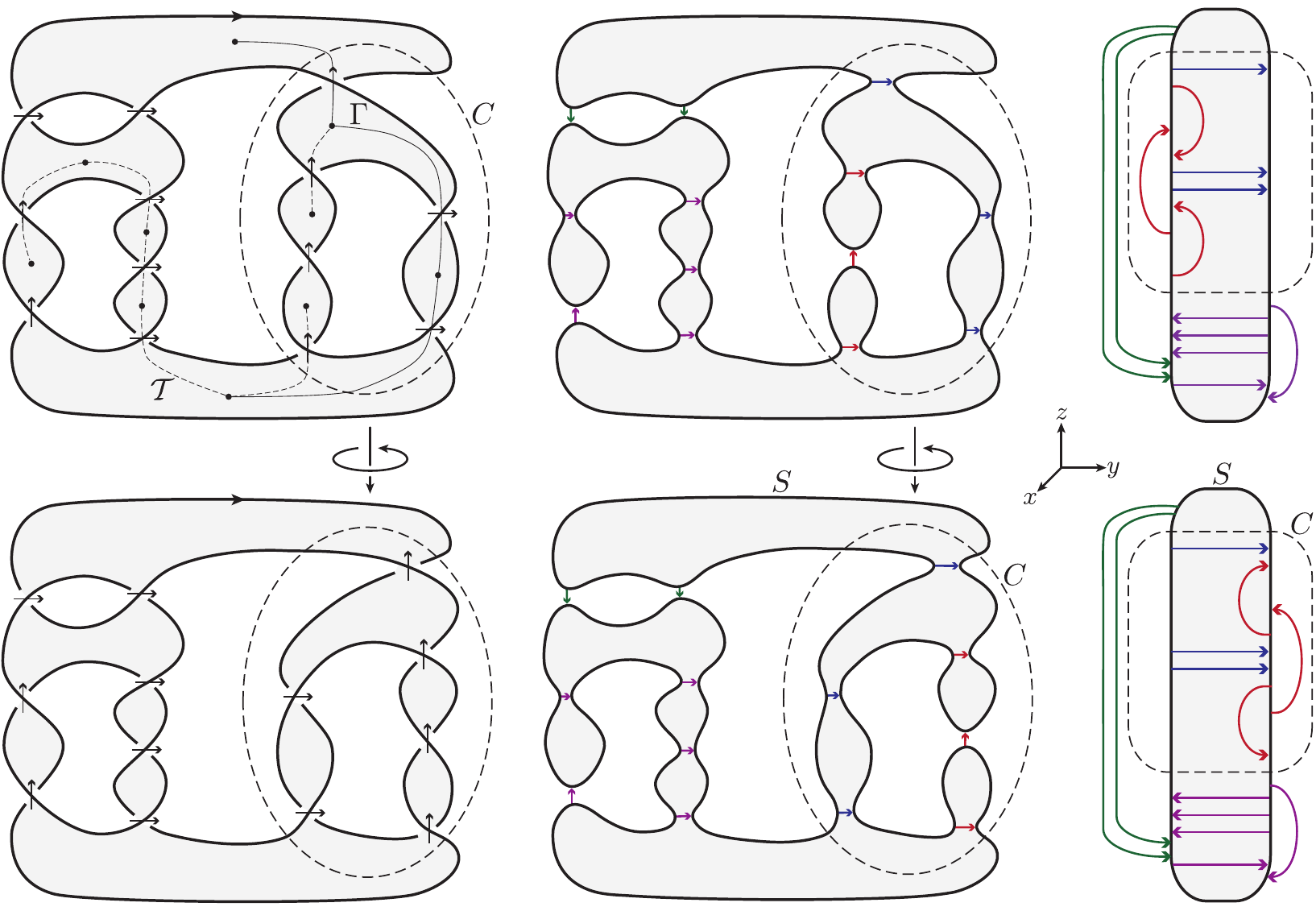}
  \caption{\textbf{The Kinoshita-Terasaka knot and the Conway knot.}  Orientations on the arcs in the upper-right resolution induce orientations on the arcs in the lower-right resolution via the mutation.  In order to obtain the same signed linking data, we have reversed the five arcs in the lower-right resolution that lie inside both $S$ and $C$.  We then work backwards to select arrows in the lower-left knot diagram.}
\label{fig:mutation}
\end{figure}

\begin{rem}
\label{rem:converse}
We have seen that if links $L$ and $L^\prime$ are related by mutation, then there exist diagrams $\dd$ and $\dd^\prime$ with crossings identified and a vertex $I^*$ such that $\dd(I^*)$ and $\dd^\prime(I^*)$ have the same mod 2 linking matrix.  The converse follows from a theorem of Chmutov and Lando \cite{chm}: \textit{Chord diagrams have the same intersection graph if and only if they are related by mutation.}  Here we view a resolution as a bipartite chord diagram, so that its mod 2 linking matrix is precisely the adjacency matrix of the corresponding intersection graph.  Note that in the bipartite case, any combinatorial mutation as defined in \cite{chm} can be realized by a finite sequence of our geometric ones.

Chmutov and Lando apply their result to the chord-diagram construction of finite type invariants.  It is known that finite type invariants of order $\leq 10$ are insensitive to Conway mutation, whereas there exists an invariant of order 11 that distinguishes the knots in Figure \ref{fig:mutation} and one of order 7 that detects genus 2 mutants (see \cite{mor} and the references in \cite{chm}).
\end{rem}

\section{The original construction of reduced odd Khovanov homology}
\label{sec:oddkh}

We now recall the original construction of reduced odd Khovanov homology, following \cite{orsz}.  Given an oriented link $L \subset S^3$, we fix a decorated, oriented diagram $\dd$ as before, though now it need not be connected.  For each vertex $I \in \vv(\dd)$, the resolution $\dd(I)$ consists of a set of circles $\{S^I_i\}$.  Let $V(\dd(I))$ be the free Abelian group generated by these circles.  The reduced group $\vo(\dd(I))$ is defined to be the kernel of the augmentation $\eta: V(\dd(I)) \to \mathbb{Z}$ given by $\sum a_i S^I_i \mapsto \sum a_i$.

Now let $\zx$ denote the free Abelian group on $n$ generators.  For each $I \in \vv(\dd)$, we have a boundary map $$d^I: \zx \to V(\dd(I))$$ given by $d^I x_i = S^I_j - S^I_k$, where $x_i$ is directed from $S^I_j$ to $S^I_k$ in $\dd(I)$.  Consider an edge $e \in \ee(\dd)$ from $I$ to $J$ corresponding to a increase in resolution at $c_i$.  If two circles merge as we move from $\dd(I)$ to $\dd(J)$, then the natural projection map $\{S^I_i\} \twoheadrightarrow \{S^J_i\}$ induces a morphism of exterior algebras.  Alternatively, if a circle splits into two descendants, the two reasonable inclusion maps $\{S^I_i\} \hookrightarrow \{S^J_i\}$ induce equivalent morphisms on exterior algebras after wedging by the ordered difference of the descendents in $\dd(J)$.  In other words, we have a well-defined map
$$\bar\partial^I_J : \Lambda^* \vo(\dd(I)) \to \Lambda^* \vo(\dd(J))$$
given by 
\begin{align*}
\bar\partial^I_J (v) &= \left\{
    \begin{array}{l}
        d^J x_i \wedge v \quad \text{if} \quad d^I x_i = 0 \in \vo(\dd(I)) \\
        \hspace{11.7mm} v \quad \text{if} \quad d^I x_i \neq 0 \in \vo(\dd(I)) \\ 
    \end{array} \right. 
\end{align*}
in the case of a split and a merge, respectively.

As in Section \ref{sec:reform}, we now form a group $\cdr$ over the hypercube and choose a type $X$ edge assignment to obtain a differential $\bar\partial_\epsilon: \cdr \to \cdr$.  The reduced odd Khovanov homology $\khor(L) \cong H_*(\cdr, \bar\partial_\epsilon)$ is independent of all choices and comes equipped with a bigrading that is initialized using $n_\pm$.  The unreduced version is obtained by replacing $\vo(\dd(I))$ with $V(\dd(I))$ above.

\begin{proof}[Proof of Proposition \ref{prop:id}]
Suppose that $\dd$ is connected.  Then for each $I\in\vv(\dd)$, the image of $d^I: \zx \to V(\dd(I))$ is precisely $\vo(\dd(I))$.  In fact, by Lemma \ref{lem:1} below, $d^I$ induces an isomorphism $\vt(\dd(I)) \cong \vo(\dd(I))$.  The collection of maps $d^I$ therefore induce a group isomorphism $\cdrr \cong \cdr$ which is immediately seen to be equivariant with respect to the edge maps $\tilde\partial^I_J$ and $\bar\partial^I_J$.  After fixing a common type $X$ edge assignment, the proposition follows.
\end{proof}

\begin{lem}
\label{lem:1}
The relations $r^I_i$ generate the kernel of the map $d^I: \zx \to V(\dd(I)).$
\end{lem}

\begin{proof}
To simplify notation, we assume that $m_i \neq m_i^*$ if and only if $i \leq k$, for some $1 \leq k \leq n$.  Consider the $n \times n$ matrix $M^I$ with column $i$ given by the coefficients of $(-1)^{m_i^*}r^I_i$.  Let $A^I$ be the leading $k \times k$ minor, a symmetric matrix.  We build an orientable surface $F^I$ by attaching $k$ 1-handles to the disk $D^2$ bounded by $S$ so that the cores of the handles are given by the arcs $x_1,\dots,x_k$ as they appear in $\dd(I^*)$.  We obtain a basis for $H_1(F^I)$ by extending each oriented arc to a loop using a chord through $D^2$.  The cocores of the handles are precisely $x_1,\dots,x_k$ as they appear in $\dd(I)$, so these oriented arcs form a basis for $H_1(F^I, \partial F^I)$.  With respect to these bases, the homology long exact sequence of the pair $(F^I, \partial F^I)$ includes the segment
\begin{align}
\label{eqn:hom}
H_1(F^I) \xrightarrow{A^I} H_1(F^I, \partial F^I) \xrightarrow{d^I |_{\zz \langle x_1,\dots, x_k \rangle}} {H}_0(\partial F^I) \xrightarrow{\eta} \zz \to 0.
\end{align}
Furthermore, the oriented chord in $D^2$ between the endpoints of $x_i$ for $i > k$ is represented in $H_1(F^I, \partial F^I)$ by the first $k$ entries in column $i$ of $M^I$.  We can therefore enlarge \eqref{eqn:hom} to an exact sequence
\begin{align*}
\zx \xrightarrow{M^I} \zx \xrightarrow{d^I} V(\dd(I)) \xrightarrow{\eta} \zz \to 0,
\end{align*}
which implies the lemma.
\end{proof}

We can reduce the number of generators and relations in our construction by using the smaller presentation in \eqref{eqn:hom}.  Namely, for each $I = (m_1,\dots,m_n) \in \vv(\dd)$, let $\vh(\dd(I))$ be the group generated by $\{x_i \, | \, m_i \neq m_i^*\}$ and presented by $A^I$.  By \eqref{eqn:rel}, the edge map $\tilde\partial^I_J$ at $c_i$ is replaced by
\begin{align*}
\hat\partial^I_J (u) &= \left\{
    \begin{array}{l}
        \hspace{3.5mm} x_i \wedge u \quad \text{if} \quad x_i = 0 \in \vh(\dd(I)) \ \text{ and } \ m_i = m_i^*\\
        -r^J_i \wedge u \quad \text{if} \quad x_i = 0 \in \vh(\dd(I)) \ \text{ and } \ m_i \neq m_i^* \\
        \hspace{11.3mm} u \quad \text{if} \quad x_i \neq 0 \in \vh(\dd(I)), \\ 
    \end{array} \right.
\end{align*}
where it is understood that $x_i \mapsto 0$ when $m_i \neq m_i^*$.  While the definition of $\hat\partial^I_J$ is more verbose, the presentations $A^I$ are simply the $2^n$ principal minors of a single, symmetric matrix: $\{(-1)^{m^*_i + m^*_j}a_{ij}\}$.  The resulting complex $(\cdrrr, \hat\partial_\epsilon)$ sits in between $(\cdrr, \tilde\partial_\epsilon)$ and $(\cdr, \bar\partial_\epsilon)$ and is canonically isomorphic to both.

\section{Branched double covers}
\label{sec:bdc}

\begin{figure}[b]
  \centering
  \includegraphics{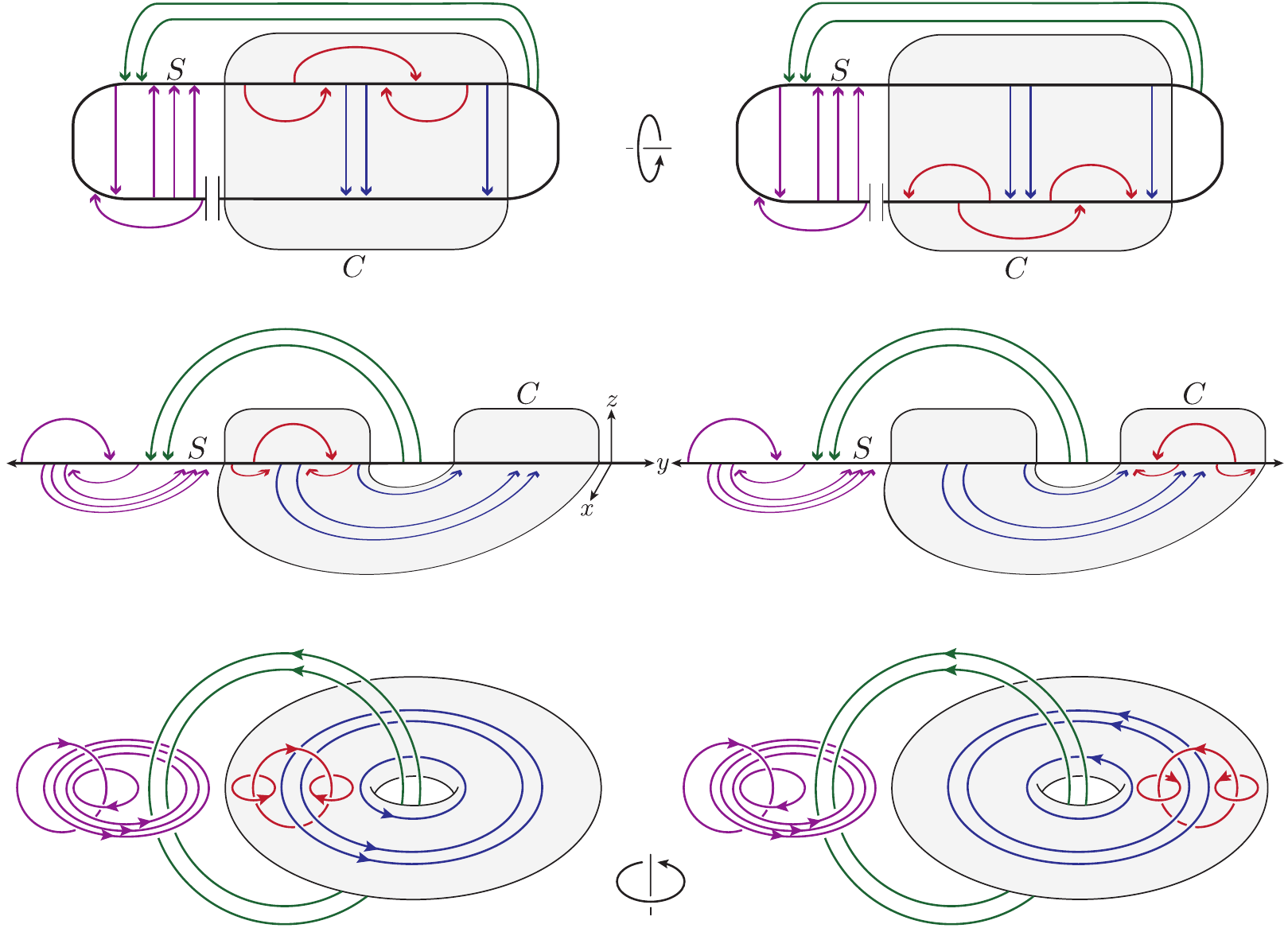}
  \caption{\textbf{Constructing a surgery diagram for the branched double cover.} The resolved diagrams in the first row are related by mutation along the Conway sphere formed by attaching disks to either side of $C$.  The double cover of $S^2$ branched over its intersection with $S$ is represented by each torus in the third row.  Rotation of the torus about the $z$-axis yields a component-preserving isotopy from $\mathbb{L}$ to $\mathbb{L}^\prime$.}
\label{fig:surgery}
\end{figure}

To a one-circle resolution $\dd(I^*)$ of a connected diagram for a link $L \subset S^3$, we associate a framed link $\mathbb{L} \subset S^3$ that presents $\Sigma(L)$ by surgery (see also \cite{gr}).  Figure \ref{fig:surgery} illustrates the construction on each of the knots in Figure \ref{fig:mutation}.  We first cut open the circle $S$ and stretch it out along the $y$-axis, dragging the arcs along for the ride.  We then slice along the Seifert surface $\{x = 0, z < 0 \}$ for $S$ and pull the resulting two copies up to the $xy$-plane as though opening a book.  This moves those arcs which started inside $S$ to the orthogonal half-plane $\{z = 0, x > 0\}$, as illustrated in the second row.  The double cover of $S^3$ branched along $S$ is obtained by rotating a copy of the half-space $\{z \geq 0\}$ by $180^\circ$ about the $y$-axis and gluing it back onto the upper half space.  The arcs $x_i$ lift to circles $\mathbb{K}_i \subset S^3$, which comprise $\mathbb{L}$.  We assign $\mathbb{K}_i$ the framing $(-1)^{m_i^*}$.

If $\dd$ is decorated, then $\mathbb{L}$ is canonically oriented by the direction of each arc in the second row of Figure \ref{fig:surgery}.  The linking matrix of $\mathbb{L}$ is then $\{-a_{ij}\}$ off the diagonal, with the diagonal encoding $I^*$.  In fact, the geometric constraints on $\mathbb{L}$ are so severe that it is determined up to isotopy by its mod 2 linking matrix (i.e., its intersection graph).  This follows from hanging $\mathbb{L}$ on a wall, or:

\begin{prop}
\label{prop:surg}
The isotopy type of $\mathbb{L} \subset S^3$ is determined by the intersection graph of $\dd(I^*)$, whereas the framing of $\mathbb{L}$ is determined by $I^*$.
\end{prop}
\begin{proof}
Suppose that $\dd(I^*)$ and $\dd^\prime(J^*)$, thought of as bipartite chord diagrams, have the same intersection graph.  Then by \cite{chm}, $\dd(I^*)$ is connected to $\dd^\prime(J^*)$ by a sequence of mutations (see Remark).   Each mutation corresponds to a component-preserving isotopy of $\mathbb{L}$ modeled on a half-integer translation of a torus $\mathbb{R}^2/\zz^2$ embedded in $S^3$ (see the caption of Figure \ref{fig:surgery}).  Therefore, the associated links $\mathbb{L}$ and $\mathbb{L}^\prime$ are isotopic.  The second statement is true by definition.
\end{proof}
\begin{proof}[Proof of Theorem \ref{thm:2}]
From the construction of the spectral sequence in \cite{osz12}, it is clear that the filtered complex associated to a connected diagram $\dd$ depends only on the framed isotopy type of the link $\mathbb{L}$ associated to some (any) one-circle resolution.  But this is preserved by mutation.
\end{proof}
\noindent Indeed,  Proposition \ref{prop:surg} implies that the link surgery spectral sequence associated to a connected diagram $\dd$ is fully determined by the mod 2 linking matrix of any one-circle resolution of $\dd$, and the vertex of this resolution (compare with Proposition \ref{prop:thrift})!

We finally come to the perspective that first motivated the construction of $(\cdrr, \tilde\partial_\epsilon)$.  For each $I = (m_1,\dots,m_n) \in \vv(\dd)$, the framing
\begin{align*}
\lambda_i &= \left\{
    \begin{array}{l}
        \infty \quad \text{if } m_i = m_i^*\\
        \hspace{1mm} 0 \hspace{.9mm} \quad \text{if } m_i \neq m_i^*
    \end{array} \right.
\end{align*}
on $\mathbb{K}_i$ gives a surgery diagram for $\Sigma(\dd(I)) \cong \#^k S^1 \times S^2$, where $\dd(I)$ consists of $k+1$ circles.  The linking matrix of $\mathbb{L}$ then presents $H_1(\Sigma(\dd(I)))$ with respect to fixed meridians $\{x_i \, | \, m_i \neq m_i^*\}$.  By identifying $H_1(\Sigma(\dd(I)))$ with $\widehat{V}(\dd(I))$, we may construct $(\cdrrr, \hat\partial_\epsilon)$, and therefore $\kho(L)$, completely on the level of branched double covers.  We will return to this perspective in \cite{jmb1}. \\

\bibliographystyle{hplain}

\bibliography{master}

\begin{thebibliography}{10}

\bibitem{b3}
J.~Baldwin.
\newblock Kh.
\newblock \url{http://www.math.princeton.edu/~baldwinj/trans.html}.

\bibitem{b}
J.~Baldwin.
\newblock On the spectral sequence from {K}hovanov homology to {H}eegaard
  {F}loer homology.
\newblock 2008, arXiv:0809.3293.

\bibitem{b2}
J.~Baldwin and O.~Plamenevskaya.
\newblock {K}hovanov homology, open books, and tight contact structures.
\newblock 2008, arXiv:0808.2336v2.

\bibitem{bar}
D.~Bar-Natan.
\newblock Mutation invariance of {K}hovanov homology.
\newblock
  \url{http://katlas.org/drorbn/index.php?title=Mutation_Invariance_of_Khovano%
v_Homology}.

\bibitem{jmb1}
J.~Bloom.
\newblock A link surgeries spectral sequence for monopole {F}loer homology.
\newblock Work in progress.

\bibitem{chm}
S.~V. Chmutov and S.~K. Lando.
\newblock Mutant knots and intersection graphs.
\newblock 2007, arXiv:0704.1313v1.

\bibitem{dun}
N.~Dunfield, S.~Garoufalidis, A.~Shumakovitch, and M.~Thistlethwaite.
\newblock Behavior of knot invariants under genus 2 mutation.
\newblock 2008, math.GT/0607258.

\bibitem{gr}
J.~Greene.
\newblock A spanning tree model for the {H}eegaard {F}loer homology of a
  branched double-cover.
\newblock 2008, arxiv/0805.1381v1.

\bibitem{kh1}
M.~Khovanov.
\newblock A categorification of the {J}ones polynomial.
\newblock {\em Duke Math. J.}, 101(3):359--426, 2000, math.QA/9908171.

\bibitem{kof}
I.~Kofman and A.~Champanerkar.
\newblock On mutation and {K}hovanov homology.
\newblock 2008, math.GT/0801.4937.

\bibitem{mor}
H.~R. Morton and N.~Ryder.
\newblock Invariants of genus 2 mutants.
\newblock 2007, arXiv:0708.0514v1.

\bibitem{orsz}
P.~Ozsv{\'a}th, J.~Rasmussen, and Z.~Szab{\'o}.
\newblock Odd {K}hovanov homology.
\newblock 2007, math.QA/0710.4300.

\bibitem{osz12}
P.~Ozsv{\'a}th and Z.~Szab{\'o}.
\newblock On the {H}eegaard {F}loer homology of branched double-covers.
\newblock {\em Adv. Math.}, 194(1):1--33, 2005, math.SG/0309170.

\bibitem{w}
L.~Watson.
\newblock A remark on {K}hovanov homology and two-fold branched covers.
\newblock 2008, arXiv:0808.2797.

\bibitem{we}
S.~Wehrli.
\newblock Khovanov homology and {C}onway mutation.
\newblock 2003, math.GT/0409328.

\end{thebibliography}

\end{document}